\documentclass[10pt,a4paper]{amsart}

\usepackage{etex}

\pdfoutput=1

\usepackage{amsmath,amstext,amssymb,mathrsfs,amscd,amsthm}
\usepackage{amsfonts}
\usepackage[all]{xy}
\usepackage{booktabs}
\usepackage{verbatim}
\usepackage{enumerate}

\usepackage{graphicx}                   
\usepackage{subfig}
\usepackage{xspace}

\newtheorem*{theorem}{Theorem}
\newtheorem{corollary}{Corollary}

\newtheorem*{predf}{Definition} 

\newtheorem*{preremark}{Remarks}  

\newtheorem*{preremark0}{Remark}  

\newtheorem*{prenotation}{Notation}

\numberwithin{equation}{section}

\newcommand\mnote[1]{}

\newcommand\mmnote[1]{}

\newcommand{\bR}{\mathbb{R}}

\newcommand\lra{\longrightarrow}

\newcommand{\map}{\mathrm{map}}

\title{A note of hyperspaces and the compact-open topology}
\author{Federico Cantero}
\date{}
\thanks{The author was funded by Michael Weiss' Humboldt professor grant. He was partially supported by the Spanish Ministry of Economy and Competitiveness under grants MTM2010-15831, MTM2013-42178-P and by the Generalitat de Catalunya as member of the team 2009 SGR 119.}
\email{\texttt{fcant{\_}01@uni-muenster.de}}
\address{{\normalfont Mathematisches Institut, Universit{\"a}t M{\"u}nster, Einsteinstra{\ss}e 62, 48149 M{\"u}nster, Germany}}

\begin{document}
\begin{abstract} We prove that the inclusion $\map(X,Y)\to \map({\mathbf{K}}(X),{\mathbf{K}}(Y))$ is continuous, where ${\mathbf{K}}(X)$ is the space of non-empty compact subsets of $X$ (also known as the hyperspace of compact subsets of $X$), and both spaces of maps are endowed with the compact-open topology.
\end{abstract}
\maketitle

Let ${\mathbf{K}}(X)$ be the space of non-empty compact subsets of $X$ endowed with the Vietoris topology, which in the case when $X$ is a metric space, it is induced by the Hausdorff distance. A subbasis of the Vietoris topology is given by the following subsets:
 \[\langle U_i\rangle = \{K\in {\mathbf{K}}(X)\mid K\subset \bigcup_i U_i, K\cap U_i\neq\emptyset~\forall i\}\]
 where $\{U_i\}$ is a finite collection of open subsets of $X$.
 
 Let $\map(X,Y)$ be the set of all continuous maps from $X$ to $Y$, endowed with the compact-open topology. A subbasis for the compact-open topology in $\map(X,Y)$ is given by the subsets
 \[W(K:U)= \{f\colon X\to Y\mid f(K)\subset U\}\]
 where $K\subset X$ is compact and $U\subset Y$ is open (or merely an element of a subbasis of the topology of $Y$).
 
There is a map
\[\phi\colon \map(X,Y)\lra \map({\mathbf{K}}(X),{\mathbf{K}}(Y))\]
that sends a map $f$ to the map $\phi(f)(K)=f(K)$.

In \cite{Mizokami}, it was proven that $\phi$ is an embedding with closed image when $\map({\mathbf{K}}(X),{\mathbf{K}}(Y))$ is endowed with the pointwise convergence topology. This, together with the fact that no article citing Mizokami's result proves continuity of $\phi$ when both spaces are endowed with the compact-open topology might suggest that $\phi$ is not continuous. This is not the case, and therefore the reason of filling this gap in the literature:
\begin{theorem} The map $\phi$ is continuous if $X$ is regular.\mnote{fc:The regularity axiom is used to simplify the basis of the compact open topology and in the proof of the assertion}
\end{theorem}
From Mizokami's result it follows immediately that
\begin{corollary} The map $\phi$ is an embedding with closed image.
\end{corollary}
One advantage of the compact-open topology is that it immediately follows that:
\begin{corollary} The evaluation map
\[\map(X,Y)\times {\mathbf{K}}(X)\lra {\mathbf{K}}(Y)\]
is continuous.
\end{corollary}
This is a corollary of the theorem above and the fact that the evaluation map $\map({\mathbf{K}}(X),{\mathbf{K}}(Y))\times {\mathbf{K}}(X)\to {\mathbf{K}}(Y)$ is continuous when the mapping space has the compact-open topology (and this map is not continuous if the mapping space has a topology coarser than the compact-open topology). The author is interested in this corollary as part of the proof in \cite{Cantero:merging} that the space of embeddings between open subsets of $\bR^n$ continuously transforms closed submanifolds into closed submanifolds, when the space of closed submanifolds is equipped with the Vietoris topology.

The theorem above can be phrased in categorical terms by saying that taking hyperspaces of compact subsets defines an \emph{topological} endofunctor of the topological category of topological spaces.

\begin{proof}
It suffices to prove that for each $f\in \mathrm{map}(X,Y)$ and for each subbasic neighbourhood $Q$ of $\phi(f)$, there are subbasic neighbourhoods $P_i$ of $f$ such that $\phi(\cap_iP_i) \subset Q $.

%

By \cite[Theorem~2.5]{Michael}, a compact subset ${\mathcal K}$ of ${\mathbf{K}}(X)$ satisfies that the union $\bigcup_{K\in {\mathcal K}} K$ is compact as well.

Therefore a subbasic open subset of $\mathrm{map}({\mathbf{K}}(X),{\mathbf{K}}(Y))$ is given by pairs $W({\mathcal K}:\langle U_i\rangle)$ with
\begin{align*}
W({\mathcal K}:\langle U_i\rangle) &= \{F:{\mathbf{K}}(X)\to {\mathbf{K}}(Y)\mid F({\mathcal{K}})\subset \langle U_i\rangle\} \\
 &=\{F\mid F(K)\subset \bigcup U_i, F(K)\cap U_i\neq \emptyset, \forall K\in {\mathcal K}, \forall i\} 
\end{align*}
Let us take such subbasic neighbourhood of $\phi(f)$, and observe that the inverse image of this neighbourhood takes the form
\begin{align*}
\phi^{-1}W({\mathcal K}:\langle U_i\rangle)&= \{g\mid g(\bigcup_{K\in {\mathcal K}} K)\subset \bigcup U_i, g(K)\cap U_i\neq \emptyset, \forall K\in {\mathcal K}, \forall i\} \\
&= \{g\mid g(\bigcup_{K\in {\mathcal K}} K)\subset \bigcup U_i\}\cap \{g\mid g(K)\cap U_i\neq \emptyset, \forall K\in {\mathcal K}, \forall i\} \\
&= W\left(\bigcup_{K\in {\mathcal K}} K:\bigcup U_i\right)\cap\bigcap_i \left\{g\mid g(K)\cap U_i\neq \emptyset, \forall K\in {\mathcal K}\right\}
\end{align*}

\noindent\emph{\underline{Assertion}}: For each $U'_i:= f^{-1}(U_i)$ there is an open subset $V:=V_i$ of $X$ such that its closure $\overline{V}_i$ is contained in $U_i'$ and for each $K\in {\mathcal K}$ it holds that $K\cap V_i\neq \emptyset$. 
\begin{itemize}
\item[]If this were not true, for each such $V$ there would exist a $K_{V}\in {\mathcal K}$ such that $K_V\cap V=\emptyset$. The open subsets of $U$ whose closure is contained in $U$ are partially ordered by inclusion, hence the assignment $V\mapsto K_V$ defines a net in ${\mathcal K}$ which has a convergent cofinal subnet because $\mathcal K$ is compact. The limit point $K_0\in {\mathcal K}$ of this subnet satisfies that $K_0\cap V = \emptyset$ for all $V\subset \overline{V}\subset U'_i$, and since $X$ is regular, the union of the $V$'s is equal to $U$. As a consequence, $K_0\cap U_i'=\emptyset$, but this contradicts the fact that $f(K_0)\cap U_i\neq \emptyset$.
\end{itemize}
Let us choose $V_i\subset U_i'$ with the above property. Take now a point $p_K^i\in V_i\cap K$ for each $i$ and $K$, and define $K_i$ to be the closure of $\{p_K^i\}$, which is compact because it is a closed subset of $\bigcup_{K\in {\mathcal K}} K$ (which is itself compact). By construction it holds that $f(K_i)\subset \bar{V}_i\subset U_i$, so $f\in W(K_i: U_i)$ and also
\[\phi(W(K_i: U_i))\subset \left\{g\mid g(K)\cap U_i\neq \emptyset, \forall K\in {\mathcal K}\right\}.\] 
As a consequece, the neighbourhood
$$
W(\bigcup_{K\in {\mathcal K}} K: \bigcup_i U_i)\cap \bigcap_i W(K_i:U_i)
$$
of $f$ is mapped into the neighbourhood $W({\mathcal K}:\langle U_i\rangle)$, and so $\phi$ is continuous.
\end{proof}

\bibliographystyle{apalike}
\bibliography{biblio-article}

\def\cprime{$'$}
\begin{thebibliography}{}

\bibitem[Cantero, 2014]{Cantero:merging}
Cantero, F. (2014).
\newblock The space of merging submanifolds in $\mathbb{R}^n$.
\newblock in preparation.

\bibitem[Michael, 1951]{Michael}
Michael, E. (1951).
\newblock Topologies on spaces of subsets.
\newblock {\em Trans. Amer. Math. Soc.}, 71:152--182.

\bibitem[Mizokami, 1998]{Mizokami}
Mizokami, T. (1998).
\newblock The embedding of a mapping space with compact open topology.
\newblock {\em Topology Appl.}, 82(1-3):355--358.
\newblock Special volume in memory of Kiiti Morita.

\end{thebibliography}

\end{document}